\documentclass[10pt]{amsart} 
\usepackage[]{amsmath, amsthm, amsfonts, amssymb}

\newcommand {\C} {{\mathbb C}}
\newcommand {\R} {{\mathbb R}}
\newcommand {\Z} {{\mathbb Z}}
\newcommand {\Q} {{\mathbb Q}}
\newcommand {\HH} {{\mathbb H}}

\newcommand {\LL} {{\mathbb L}}

\newcommand {\E} {{\mathcal E}}

\newcommand {\G} {{\mathcal G}}
\newcommand {\B} {{\mathcal B}}
\newcommand {\F} {{\mathcal F}}
\newcommand {\del} {\nabla}

\newcommand {\dt} {\bullet}
\newcommand {\bF} {{\widehat{F}}}
\newcommand {\V}  {{\mathcal V}}

\newcommand {\A}[1] {A_X^{#1}(\log\, D)}

\newcommand {\ju}[1] {j_{U_{#1}*}^*}
\newcommand {\Cf} {{\mathcal C}{\mathcal F}}
\newcommand {\HC} {{\mathbb H}{\mathbb C}}
\newcommand {\SV} {{\mathcal S}{\mathcal V}}

\newtheorem{thm}[subsection]{Theorem}

\newtheorem{lemma}[subsection]{Lemma}
\newtheorem{prop}[subsection]{Proposition}

\newtheorem{remark}[subsection]{Remark}
\newtheorem{ex}[subsection]{Example}

\title{ Building Mixed Hodge Structures}
\author{Donu Arapura}
\address{Department of Mathematics, Purdue University, West Lafayette
  IN 47907}
\thanks{Partially supported by NSF}

\begin{document}
\maketitle

This article is aimed at people who already know what mixed Hodge
structures are and what they are good for, but who are not sure how 
to construct them. 
Even if one wants to treat mixed Hodge theory primarily as a black box, 
there will be times when one needs to open it up, look inside, and 
perhaps reassemble it in a different way.  I would like 
to discuss 
some of the  mechanics involved along the lines of  Deligne's original
construction \cite{D}. I should mention that a somewhat different
and far reaching approach to constructing mixed Hodge structures has
been given by M. Saito \cite{saito}; I will say a few
words about it at the end.  However, I believe that there is still some 
value  in this old-fashioned approach. Some of the techniques
involved are quite useful in other situations, 
especially descent theory,
which has  a number of  applications
beyond those discussed here (e.g. \cite{beilinson},
\cite{gillet-soule}, \cite{gnpp}, \cite{guillen-navarro},
\cite{hain}, \cite{navarro}).

The first basic principle  of the construction
is that in order to get mixed Hodge structures
on cohomology
of a space, 
the necessary filtrations should  be defined on the
level of  complexes
that give rise to the cohomology. 
Appropriate homological conditions need to be imposed in 
order to get things to  work. In this article, these conditions will be 
embodied in  the notion of a $\C$-Hodge complex  which is a fairly 
minimalistic reworking of Deligne's original axioms for a cohomological 
mixed Hodge complex.  $\C$-Hodge complexes need not have any underlying 
integral or real structure, but they do possess the necessary structure 
to carry out the key degeneration of spectral sequence arguments. Once these
technical results are established, it is fairly 
straightforward 
to impose integral structures on the $\C$-Hodge complexes 
so as to obtain mixed Hodge structures on cohomology.

The second basic principle is that in order to construct mixed Hodge 
structures on the cohomology of a  whole class of spaces and
coefficient sheaves, it is only necessary to construct 
explicit Hodge complexes on a  subset of these pairs, provided that the
general pairs can be 
``approximated'' 
by the special ones.
This reduction is ultimately  handled by cohomological descent.
However, in order to clarify the structure of the
argument, the condition  for guaranteeing the existence
of approximations, called density, will be introduced first;
the actual descent theory will be done afterwards.
Verification of density in practice requires resolution of 
singularities. 

The last section contains a rather elaborate illustration of these 
ideas in construction of mixed Hodge structures on the cohomology
of certain constructible sheaves.

My heartfelt thanks goes out to the referees for their remarkable
thoroughness. I have learned after the fact, that one of them was
Steve Zucker.


\section{First attempt}

Let me begin by trying to make a bare hands construction of a
functorial mixed Hodge structure on the cohomology of a 
singular projective variety. This attempt  will fail, but
it will be instructive 
nevertheless.

For the record, I should recall what a mixed Hodge structure is.
I will use the 
standard
convention that filtrations indexed by superscripts
are decreasing while 
those with
subscripts are increasing. These can be
interchanged by the rule $W^{-n} = W_n$,
and we do this without further comment.
A mixed Hodge structure is usually taken to be a finitely generated
Abelian group $H_\Z$, along with an ascending filtration 
$W_\dt$ on $H_\Q= H_\Z\otimes \Q$, and a decreasing filtration 
$F^\dt$ on $H=H_\Z\otimes \C$ subject to certain axioms \cite{D}. After
jettisoning $H_\Z$ and $H_\Q$, one arrives at the notion
of the {\em $\C$-mixed Hodge structure} (or simply just a 
$\C$-Hodge structure) which is a finite dimensional vector space $H$ together
with  three filtrations $(F^\dt,  \bF^\dt, W_\dt)$
 which  are opposed in the
sense of Deligne ~\cite{D}, i.e., 
$$ Gr^W_n H = F^pGr^W_nH \oplus \bF^qGr^W_nH$$
whenever $p + q = n$,
where $Gr^W_n H = W_n/W_{n-1}$.
By tradition  $W$ is taken to be increasing, but it will be frequently
convenient to regard a $\C$-Hodge structure as given by three
decreasing filtrations $(F^\dt,  \bF^\dt, W^\dt)$.
 The third filtration $\bF^\dt$ plays the role of the
conjugate filtration to $F^\dt$.  Morphisms of $\C$-Hodge structures
are  simply
linear maps preserving the filtrations.  
In the presence of several $\C$-Hodge structures $H, H'\ldots$,
I will indicate the corresponding filtrations by 
$F^\dt H, F^\dt H',\ldots $. 
$\C$-Hodge structures possess most 
of the homological properties of mixed Hodge structures. For example the 
category of $\C$-Hodge structures is Abelian, and the functors 
$Gr_F,Gr_{\hat F}, Gr_W$ are exact \cite[1.2.5]{D}.
A $\C$-Hodge structure is pure of weight $n$ if $Gr_n^WH = H$.
The filtration $W$ will usually be suppressed for pure Hodge
structures.  There are no nonzero morphisms  between pure $\C$-Hodge structures of
different weights.  An integral lattice on a $\C$-Hodge
structure $(H,  F^\dt, \bF^\dt, W_\dt)$
 is a finitely generated Abelian group
$H_\Z$ with an isomorphism  $H_\Z\otimes \C \cong H$
such that $W_\dt$ is rational, 
that is spanned by $W_\dt \cap (H_\Z\otimes \Q)$,
and  $\bF^\dt= {\overline F}^\dt$ the conjugate of $F$ with respect to 
the induced real structure on $H$. A mixed Hodge structure 
in the usual sense is the same 
thing as a $\C$-Hodge structure together with a lattice. A
morphism of mixed Hodge structures is simply a morphism of $\C$-Hodge
structures and a compatible map of lattices.

I will generally use the term complex variety to mean
a reduced separated scheme of finite type over $\C$
with its classical topology.  Suppose that $X$ is complex  
proper variety. By Hironaka's theorem \cite{hironaka}, there exists 
a commutative diagram
$$
\begin{array}{ccc}
E& \stackrel{i}{\hookrightarrow} & \tilde X\\
\downarrow p &        & \downarrow \pi \\
S& \stackrel{j}{\hookrightarrow} & X\\
\end{array}
$$
where $\pi$ is a resolution of singularities with exceptional
divisor $E$, and $S = \pi(E)$. Then

\begin{prop}\label{prop:mayervit}
        For any sheaf of Abelian groups $\F$ on $X$,
        the ``Mayer-Vietoris'' sequence 
$$\ldots H^i(X,\F) \stackrel{\pi^*+i^*}{\longrightarrow}
 H^i(\tilde X, \pi^*F)\oplus H^i(S,j^*\F)
\stackrel{i^*-p^*}{\longrightarrow}
 H^i(E, (\pi\circ i)^*\F)\to \ldots $$
is exact.
\end{prop}

Although 
this fact
is well known, I will give the proof since
it foreshadows the results of section \ref{sect:desc}.

\begin{proof}
I will use some basic notions from derived categories \cite{verdier}
for succinctness, but the argument can  be rewritten 
without them.

Let $g= \pi\circ i$.
It is enough to show that the sequence
$$\F \to i_*\F\oplus \R \pi_*\F \to \R g_*\F$$
is part of a distinguished triangle. 
This can be checked  at the stalks. If $x\notin S$, then
the stalk of the sequence at $x$ reduces to 
$$\F_x \stackrel{=}{\to} 0\oplus \F_x\to 0$$
which is certainly part of a distinguished triangle.
If $x\in S$, let $E_x$ denote
the fiber  $g^{-1}(x) = \pi^{-1}(x)$. Then the stalk
of the sequence decomposes as sum of sequences
$$\F_x\stackrel{=}{\to} \F_x \to 0$$
and
$$0\to \R\Gamma (E_x,\F|_{E_x}) \stackrel{=}{\to} 
\R\Gamma (E_x,\F|_{E_x}) $$
which both extend to distinguished triangles.

\end{proof}

The proposition suggests a plan for putting a mixed Hodge
structure on $H^*(X,\Z)$.
First use classical Hodge theory to put  pure Hodge structures
on the cohomology of smooth projective varieties. Then
assume by induction that the cohomology of all proper varieties of
dimension less than $dim \,X$ has been equipped with a mixed
Hodge structure.
Choose a diagram as above, and try to impose a mixed Hodge structure
on it so that the Mayer-Vietoris sequence is an exact sequence
of mixed Hodge structures.
Well that's the idea.
The problem is that the category of mixed Hodge structures has
nontrivial extensions, so compatibility with Mayer-Vietoris
will not uniquely determine the mixed Hodge structure on $H^*(X)$.
Even if one  worked in the category of $\C$-mixed
Hodge structures, which is semisimple, one would only be able to
determine the $\C$-mixed Hodge structure on $H^*(X)$ up to
noncanonical isomorphism.
To circumvent these problems, it will be necessary to work
on the level of complexes. Once this is done,
it is  possible to fix up this
 inductive approach. However it
will turn out to be simpler to do the construction in one step using
descent theory (following Deligne's idea). An elementary construction
along  lines similar to the above
sketch has been given by El Zein \cite{E}. A rather
different construction has been given by Saito \cite{saito1, saito}.

\section{$\C$-Hodge Complexes}

Recall that a filtration on a complex $A^\dt$ 
is a sequence of filtrations $F^\dt A^n \subset A^n$
compatible with the differentials.
The filtration is biregular if for
each $n$, there exist $a$ and $b$ such that $F^a A^n = A^n$
and $F^b A^n = 0$. In the sequel, I will try to conserve notation
by writing $F^\dt$ or simply $F$ instead of $F^\dt A^\dt$
when no confusion is likely.

The basic device for constructing $\C$-Hodge structures is a
 {\em $\C$-Hodge complex} \cite{A} 
 It is a biregularly trifiltered complex of 
$\C$-vector spaces $(A^\dt,F^\dt A^\dt  , \bF^\dt A^\dt , W^\dt A^\dt)$ 
bounded below and
satisfying: 
 
\begin{enumerate} 

\item[HC1.] $A^\dt$ has finite dimensional cohomology.

\item[HC2.] For each $p$, the filtrations on $Gr^p_WA$ induced by $F^\dt$
and $\bF^\dt$ are
strictly preserved by the differentials (i.e. 
$dF^p = image(d)\cap F^p$ for all $p$, and likewise for $\bF^\dt$).

\item[HC3.] $H^i(Gr^p_W A)$, when equipped with the filtrations
$$image[H^i(F^\dt Gr^p_W A)\to H^i(Gr^p_W A)]$$
$$image[H^i(\bF^\dt Gr^p_W A)\to H^i(Gr^p_W A)],$$
becomes a pure $\C$-Hodge structure of weight $i-p$ 
for all $i,p$.

\end{enumerate}

This is related to 
Deligne's notion of a mixed Hodge complex \cite{D}; a more
detailed comparison will be given below.

\begin{ex}\label{ex:easy}
Consider a complex $A^\dt$ of $\C$-mixed Hodge structures with
its filtrations $F^\dt, \bF^\dt, W^\dt$. As it stands, this is
not a $\C$-Hodge complex. But if one defines 
$$\tilde W^iA^n = W^{i-n}A^n$$
then $(A^\dt,F^\dt, \bF^\dt, \tilde W^\dt )$ is.

\end{ex}

\begin{ex}\label{ex:trans}
Given a $\C$-Hodge complex $(A^\dt, F^\dt , \bF^\dt, W^\dt)$,    
let $(A[i]^\dt, d_{A[i]}) = (A^{\dt+i}, (-1)^id_A)$
be the translated complex. Then this becomes a $\C$-Hodge complex
 when endowed with filtrations $F$, $\bF$ and $W^\dt[i] = W^{\dt+i}$.
For any pair of integers the twist $A(a,b)$ which is
$A$ with the filtrations $F^\dt[a]$, $\bF^\dt[b]$ and $W^\dt[-a-b]$ is
a $\C$-Hodge complex. $A(1,1)$ is the Tate twist.
\end{ex}

\begin{ex}\label{ex:kahler} Let $X$ be a compact K\"ahler manifold. 
The de Rham complex $\E^\dt(X)$ with
filtrations
$$F^p\E^\dt(X) = \bigoplus_{p'\ge p}\, \E^{p'q}(X)$$
$$\bF^q\E^\dt(X) = \bigoplus_{q'\ge q}\, \E^{pq'}(X)$$
$$W_0\E^\dt(X) = \E^\dt(X);\> W_{-1}\E^\dt(X) = 0$$
is a $\C$-Hodge complex by classical Hodge theory. 
\end{ex}

\begin{ex}\label{ex:vhs}
Let $X$ be as above.
A polarized complex variation of Hodge structure is a $C^\infty$ vector 
bundle $V$ with
an indefinite $C^\infty$ Hermitian metric $<,>$, a decomposition
$V=\bigoplus_{p+q=w}\, V^{pq}$,  and a flat connection $\nabla$ on $V$
subject to the following axioms:
\begin{itemize}

\item The  metric  is positive definite on $V^{pq}$ for $p$ even and
negative definite when $p$ is odd, and the decomposition
is orthogonal. 

\item The connection is compatible with the metric (i.e. $\nabla$
satisfies the Leibniz rule with respect to $<,>$) and satisfies
the following form of  Griffith's transversality:
$$\nabla(V^{pq}) \subseteq \E^{01}(V^{p+1,q-1}) \oplus
\E^{10}(V^{p,q}) \oplus \E^{01}(V^{p,q}) \oplus
\E^{10}(V^{p-1,q+1})$$
where $\E^{rs}(U)$ is the sheaf $C^\infty$ $U$-valued $rs$-forms.

\end{itemize}

The de Rham complex $(\E^\dt(X,V),\nabla)$ becomes a $\C$-Hodge
complex when equipped with the filtrations
$$F^p\E^\dt(X,V) = \bigoplus_{p'+p''\ge p}\, \E^{p'q'}(X,V^{p''q''})$$
$$\bF^q\E^\dt(X) = \bigoplus_{q'+q''\ge q}\, \E^{p'q'}(X,V^{p''q''})$$
$$W_w\E^\dt(X) = \E^\dt(X,V);\> W_{w-1}\E^\dt(X) = 0.$$
This follows from the generalized K\"ahler identities of Deligne and
Simpson \cite{simpson}, \cite{zucker}.
A lattice is  a locally constant sheaf ${\mathcal V}$ of Abelian
groups with finitely generated
stalks  and an isomorphism ${\mathcal V}\otimes \C \cong ker(\nabla)$
such that $V^{pq}$ is conjugate to $V^{qp}$ with respect to resulting
$\R$-structure. The usual notion of a polarized variation
of Hodge structure is equivalent to a complex variation of Hodge
structure and a lattice.
See \cite{D2} \cite{griffiths-schmid}, \cite{simpson} for more information.

\end{ex}

 Complex (usually nonintegral) variations of Hodge structure,
 occur naturally in
Simpson's work \cite{simpson}
as fixed points in the moduli space of Higgs bundles
under the natural $\C^*$-action. Using this, he has shown that
any representation of the fundamental group of
a smooth projective variety can be deformed
to the monodromy of a complex variation of Hodge structure.

In the previous examples, the complexes arose by taking global
sections of complexes of sheaves. This can be generalized
as follows.
A  $\C$-Hodge complex on a space $X$ is a quadruple
 $(A^\dt,F^\dt, \bF^\dt, W^\dt )$
consisting of   a bounded below complex of sheaves of 
$\C$-vector spaces $A^\dt$ with finite dimensional  cohomology, together
with  three biregular filtrations on it such that:

\begin{enumerate} 

\item[$\HC 1$.] $A^\dt$ has finite dimensional hypercohomology.

\item[$\HC 2$.] For each $p$, the filtrations on ${\mathbb R}\Gamma Gr^p_WA$
  induced by $F$ and $\bF$ are strictly preserved by the differentials.

\item[$\HC 3$.] For every $i,p$, $F$ and $\bF$ induce  a pure $\C$-Hodge 
structure of 
weight $i-p$ on $\HH^i(Gr^p_W A)$. 

\end{enumerate}

A map of trifiltered complexes
$$(A^\dt,F^\dt A^\dt, \bF^\dt A^\dt, W^\dt A^\dt) \to
(B^\dt, F^\dt B^\dt, \bF^\dt B^\dt, W^\dt B^\dt)$$
will be called an {\em acyclic resolution}
if the associated
maps  
$$Gr_{FA}^\dt Gr_{WA}^\dt A^\dt\to Gr_{FB}^\dt Gr_{WB}^\dt B^\dt$$ 
and 
$$Gr_{\hat FA}^\dt Gr_{WA}^\dt A^\dt\to 
Gr_{\hat FB}^\dt Gr_{WB}^\dt B^\dt$$
are acyclic
resolutions in the usual sense. Canonical acyclic resolutions can 
 be constructed  using Godement's flasque resolutions.
The axioms imply that if $(A^\dt,F^\dt, \bF^\dt, W^\dt )$ is a 
$\C$-Hodge complex, then after replacing it by
an acyclic resolution, $\Gamma(A^\dt)$ with its induced
filtrations is
a $\C$-Hodge complex in the earlier sense.

\begin{ex}\label{ex:log}
Let $D\subset X$ be a divisor with normal crossings in a 
compact K\"ahler manifold and let $j:X-D\hookrightarrow X$ be the
inclusion. 
Let $A_{X}^\dt$ be the sheaf of $\C$-valued real analytic differential
forms on $X$. 
Following Navarro Aznar \cite{navarro}, define $A_X^\dt(log\, D)$ to be the 
$A_X^0$-subalgebra
of $j_*A_{X-D}^{\dt}$ generated by $A_X^\dt$ and $\log|f|$, $
Re\,df/f$, and  $Im\, df/f$  for
any (local) holomorphic function $f$ vanishing along a component of $D$.
This is a subcomplex of $j_*A_{X-D}^\dt$, which is an
incarnation of ${\mathbb R}j_*\C$. Furthermore
it has a bigrading 
$$A_{X}^n(logD) = \bigoplus_{p+q=n}\, A_{X}^{p,q}(logD)$$ 
by $(p,q)$ type.
 With  filtrations 
$$F^p \A {\dt} = \bigoplus_{p' \ge p} \, \A{p',q} $$ 
$$ \bF^q \A {\dt} = \bigoplus_{q' \ge q} \, \A {p,q'} $$
$$ W_k A_{X}^{m}(\log D) = image(A_{X}^{k}(\log D) \wedge A_{X}^{m-k} \to 
A_{X}^{m}(\log D) ) $$ 
$A_X(\log D)$ becomes a $\C$-Hodge complex. The complex
$A_X^\dt(\log D)$  can be replaced by the $C^\infty$ analogue 
 constructed by Burgos \cite{burgos} to obtain an explicit
acyclic resolution of this $\C$-Hodge complex.
\end{ex}

\begin{ex}\label{ex:ulog}
Let $X$ and $D$ be as above.
Let $(V^o,\del^o)$ be a holomorphic bundle on $X-D$ with a 
flat connection which is unitary in the sense
 that $\V = ker(\del^o)$ corresponds to a unitary representation
of $\pi_1(X-D)$. 
Then there exists a unique extension of $V^o$ to a vector
bundle $V$ on $X$ such that $\del^o$ extends to a logarithmic
connection $\del$ with residues having eigenvalues in $[0,1)$ 
\cite{regsing}. The  complex $(\A {\dt} \otimes V,\del)$ is
quasi-isomorphic to
$\R j_* \V$. There is a filtration $W$ on 
$(\A {\dt} \otimes V,\del)$, such that together with the filtrations
$$F^p \A {\dt}\otimes V = \bigoplus_{p' \ge p} \, \A{p',q}\otimes V  $$ 
$$ \bF^q \A {\dt}\otimes V = \bigoplus_{q' \ge q} \, \A {p,q'}\otimes V$$
this becomes a $\C$-Hodge complex \cite[IV.1]{A}.
\end{ex}

\begin{ex}\label{ex:cone} 
A morphism of $\C$-Hodge complexes is a map of complexes which 
is compatible with 
the filtrations. Let $f:(A^\dt, F^\dt A^\dt,\ldots)\to 
(B^\dt, F^\dt B^\dt,\ldots)$
be a morphism. The mapping cone is a new complex $Cone(f)^\dt = C^\dt$
with
$$C^n = A^{n+1}\oplus B^{n}$$
and $d(\alpha,\beta) = (-d\alpha, d\beta-f(\alpha))$.
 With the filtrations
$$F^p C^n = F^p A^{n+1} \oplus F^pB^{n}$$
$$\bF^p C^n = \bF^p A^{n+1} \oplus \bF^p B^{n}$$
and
$$W_k C^n = W_{k-1} A^{n+1} \oplus W_{k}B^{n}$$
$C^\dt$ becomes a $\C$-Hodge complex.
By construction there are morphisms
$B\to Cone(f)$ and $Cone(f)\to A[1]$ where $A[1]$ is defined
in example \ref{ex:trans}
\end{ex}

The previous example can be generalized.

\begin{ex}\label{ex:total}
 A differential
graded $\C$-Hodge complex consists of a collection  $(A^{ n\dt}, 
FA^{n\dt},\ldots)$
of $\C$-Hodge complexes, and morphisms
$\delta_n:A^{n, \dt}\to A^{n+1, \dt}$ satisfying 
$\delta_{n+1}\delta_n=0$.
Define the total (or associated single) complex $T^\dt=tot^\dt(A)$
by 
$$T^n = \bigoplus_{i+j=n} A^{ij}$$
with differential $d +(-1)^j\delta_j$. With filtrations
$$F^pT^n = \bigoplus_{i+j=n} F^pA^{ij}$$
$$\bF^pT^n = \bigoplus_{i+j=n} \bF^pA^{ij}$$
and
$$ W_k T^n = \bigoplus_{i+j=n} W_{k+i}A^{ij}$$
$T^\dt$ becomes a $\C$-Hodge complex.
\end{ex}

The above constructions will be used extensively later. Here
are a few concrete instances.

\begin{ex}\label{ex:jshriek}
Let $X$ be a compact K\"ahler manifold, $i:Y\hookrightarrow X$ a closed
submanifold and $j: U = X-Y\hookrightarrow X$ the 
the complement. Let $C^\dt$ be the mapping cone of the
restriction map $\E_X^\dt\to \E_Y^\dt$. This is a $\C$-Hodge
complex quasi-isomorphic to $j_!\C_U$, and $H^i(j_!\C_U)\cong H^i(X,Y)$.
\end{ex}

 The above example is very special case of the constructions
of the final section.

\begin{ex}\label{ex:theta}
Let $X$ be a compact K\"ahler manifold, and let $\theta\in 
H^0(X,\Omega_X^1)$ and $N$ be some positive integer.
Let $\theta'$ denote the operator $\phi\mapsto 
(-1)^{deg\phi}\theta\wedge \phi$. Then
 $(\theta')^2=0$, we get a differential
graded $\C$-Hodge complex
$$\E_X^\dt\stackrel{\theta'}{\longrightarrow}\E_X^\dt(1,0)[1]
\stackrel{\theta'}{\longrightarrow} \ldots\E_X^\dt(N,0)[N]$$
where $\E_X(i,0)[i]$ is defined in \ref{ex:trans} and \ref{ex:kahler}.
\end{ex}

This example is from \cite{A}.
Related constructions occur in the work of 
Steenbrink \cite{steenbrink}.

In example \ref{ex:easy}, the cohomology groups of $A^\dt$
carry $\C$-Hodge structures with filtrations induced by 
$F$, $\bF$ and $W$,
because the category of $\C$-Hodge structures is Abelian. 
This generalizes to all
$\C$-Hodge complexes. The following, which  was announced in \cite{A}
and is an analogue of \cite[scholium 8.1.9]{D}, is the basic result.

\begin{thm}\label{thm:1}
 Let $(A^\dt,F^\dt, \bF^\dt, W^\dt )$ be a $\C$-Hodge complex 
then \begin{enumerate}
        \item[(a)] The spectral sequence induced by $W$ degenerates at 
         $E_2$.
        \item[(b)] $H^i(A^\dt)$ carries a $\C$-Hodge structure with
        filtrations given by $F^\dt,\, \bF^\dt$ and
        $W^\dt$ shifted by $i$,  or explicitly by:
        $$F^pH^i(A^\dt) = image[H^i(F^pA^\dt) \to H^i(A^\dt)]$$
        $$\bF^pH^i(A^\dt) = image[H^i(\bF^pA^\dt) \to H^i(A^\dt)]$$
        $$W^{k}H^i(A^\dt) = image[H^i(W^{k+i}A^\dt) \to H^i(A^\dt)]$$

       \item[(c)] The spectral sequences associated to $F$ and $\bF$
         degenerate at $E_1$.  
\end{enumerate}
\end{thm}

The proof of the above theorem will be given in Section 3. 
A variant of part (a), proved by a similar argument, is:

\begin{prop}\label{prop:thm1}
 The spectral sequence induced by $W^\dt$ 
 on $Gr_F^iA^\dt$ degenerates   at $E_2$.
\end{prop}

Analogues of (c) also hold,
for $F$ and $\bF$,
namely  the spectral sequences induced by $F^\dt$ and $\bF^\dt$ on 
$Gr_W^i A^\dt$ degenerate at $E_1$. This is nothing but a restatement
of axiom HC1.
As a corollary, one gets similar statements for $\C$-Hodge complexes
on topological space with cohomology replaced by hypercohomology. 

The degeneration statements have a number of important applications.
When part (c) of the theorem
 is applied to example \ref{ex:ulog}, along 
with the filtered quasi-isomorphism $(\Omega_X^\dt(\log D)\otimes V, F^\dt)
\hookrightarrow (\A {\dt} \otimes V, F^\dt)$,
 one obtains the degeneration of the spectral sequence
$$H^q(\Omega_X^p(\log D)\otimes V)\Rightarrow H^{p+q}(X-D,\, \V)$$
This fact, and its generalizations,
 can be used to prove   vanishing theorems (see 
\cite{esnault-v}). For example, if $X-D$ is Stein (e.g.
if $X$ is projective and $D$ ample) then
it has the homotopy type of a CW-complex of dimension at
most $n=dim\,X$; this implies the following version of
Kodaira's vanishing theorem: 
$$H^i(X, \omega_X(D)\otimes V)\subseteq H^{n+i}(X-D, \V) = 0$$
for $i>0$. 
When proposition \ref{prop:thm1} is  applied to  the complex
of    example \ref{ex:theta}, one recovers 
a theorem of Green and Lazarsfeld \cite{green-l} that the
spectral sequence
$$E_1= H^q(X,\Omega_X^p)\Rightarrow 
\HH^{p+q}(X;O_X\stackrel{\theta}{\longrightarrow} \Omega_X^1
\stackrel{\theta}{\longrightarrow}\ldots)$$
degenerates at $E_2$. 
This kind of result, while technical,
is crucial for establishing their
generic vanishing theorems.
The argument is adapted from \cite{A} which contains further
details and refinements.

With the basic homological apparatus  in place,  the  integral
structure can now be added. 
An integral structure for a $\C$-Hodge complex
$(A^\dt, W,F,\bF)$ on a topological space  consists of
a bounded below complex of sheaves of Abelian groups $A^\dt_\Z$
with finitely generated cohomology and 
quasi-isomorphisms
$$ A_\Z^\dt\otimes\Q  \stackrel{\sim}{\leftarrow} { A'}^\dt
\stackrel{\sim}{\to} A_\C^\dt$$
inducing lattices on cohomology. Note that the intermediate
complex ${A'}^\dt$ is part of the data (or to 
put it another way: fix an explicit representative
for an isomorphism $ A_\Z^\dt\otimes\Q  \cong  A_\Q^\dt$
 in the derived category ). The easist way to ensure
that $W_\dt H^\dt(A^\dt)$ is rational
is to construct a filtered complex
of $\Q$-vector spaces $(A^\dt_\Q, W_\Q)$ and 
isomorphisms 
$$ A_\Z^\dt\otimes\Q  \cong A_\Q^\dt$$
$$ (A_\Q^\dt,W_\Q)\otimes \C \cong (A^\dt,W)$$
in the derived and filtered derived categories.
A Hodge complex
will consist of a $\C$-Hodge complex
together with a choice of integral structure. 
Hodge complexes  carry mixed Hodge structures on cohomology.
A trivial example of an integral structure is the inclusion of the
constant sheaf $\Z_X\subset \E_X^\dt$ into the De Rham complex of 
a compact K\"ahler manifold. More generally for  variations
of Hodge structure, the inclusion of the lattice in $\E_X^\dt(V)$
supplies an integral structure.
In example \ref{ex:log}, $A_\Z$ can be taken to be the
direct image under $j_*$ of the complex
of sheaves of singular  $C^\infty$ cochains,
$A' = A=A_X^\dt(log\, D)$, and $A'\to A_\Z\otimes \C$
the map which assigns to a form $\alpha$, the cochain 
$\gamma\mapsto \int_\gamma\alpha$ (the rationality of $W$
is a consequence of  a filtered quasi-isomorphism 
$(A_\Z, \tau_{\le\dt})\otimes \C \cong (A, W_\dt)$,
\cite{D, navarro}).

A few remarks ought to be made  about the relationship
between Hodge complexes 
and (cohomological) mixed Hodge complexes \cite[8.1]{D}. The set up
for the latter is similar, one has a collection of complexes with 
filtrations
and compatability isomorphisms subject to appropriate conditions.
Although I do not expect a  precise equivalence,
it does seem to be case in practice that whenever one has a natural
Hodge complex, one has a corresponding cohomological mixed Hodge
and conversely. However, there are interesting examples of 
$\C$-Hodge complexes
(or structures) which do not 
possess any integral or even real structure.
The $\C$-Hodge complexes of  examples  \ref{ex:vhs},
\ref{ex:ulog} and \ref{ex:theta} will not  generally possess 
natural integral structures.
For certain applications (such as 
those
discussed
above), the degeneration of the above  spectral
sequences is precisely what is required; the integral structure is
irrelevant even when it exists.
For these reasons, 
the notion of a $\C$-Hodge complex is useful on its own.

For any space $X$, the collection of $\C$-Hodge complexes
forms a category, 
denoted 
$C_{\C-Hodge}(X)$. More generally,
let $C_{\C-Hodge}$ be the category whose objects are pairs consisting
of a topological space and a $\C$-Hodge complex on it. A morphism from 
$(X,(A,\ldots))$
to $(Y, (B,\ldots))$ is pair consisting of a continuous map
$f:X\to Y$ and a map of trifiltered complexes $(f^*B, \ldots)\to
(A,\ldots)$. 
This category is fibered over  the category of topological spaces
by the functor $|\ldots |$,
and $C_{\C-Hodge}(X)$  is the fiber over $X$.
The category of Hodge complexes $C_{Hodge}$ is defined similarly.

\section{ Proof of Theorem \ref{thm:1}}

Let $A^\dt$ be a bounded below complex of vector spaces with two
biregular filtrations  $W^\dt, F^\dt$.
There is a spectral sequence corresponding to the first filtration with
$$E_1^{pq}(A,W)  = H^{p+q}(Gr_W^p A^\dt)$$
and differential 
$$d_1:E_1^{pq}(A,W) \to E_1^{p+1,q}(A,W)$$
given by the connecting map associated to 
$$ 0\to Gr_W^{p+1}A \to W^pA/W^{p+2}A \to Gr_W^pA\to 0.$$
This can be filtered  by $F$ in many (usually) different ways. The  
{\em first direct filtration} is
$$F_{dir}^iE_r^{pq}(A,W) = image[E_r^{pq}(F^i,F^i \cap W) 
\to E_r^{pq}(A,W)].$$
The {\em recursive filtration} is defined inductively by 
$F_{rec}^i E_1^{pq} = F_{dir}^i E_1^{pq}$ and
$$F_{rec}^i E_{r+1}^{pq}(A,W) = image[F_{rec}^i\cap ker_r^{pq}(A,W)]$$
where 
$$ ker_r^{pq}(A,W) =  ker[E_r^{pq}(A,W) \to E_r^{p+r,q-r+1}(A,W)]$$
and the image is taken in
$ E_{r+1}^{pq}(A,W)$.
$E_\infty$ is a subquotient of $H^*(A)$ and it inherits a filtration
induced by
$$F^iH^*(A) = image[H^*(F^iA) \to H^*(A)]$$
which will be called the {\em final filtration}.
The first direct filtration is useful because it is preserved by 
differentials
$d_r:E_r^{pq}(A,W) \to E_r^{p+r,q-r+1}(A,W)$,
while the recursive filtration is 
what one is seeing if one follows the spectral sequence from 
$E_1$ to $E_2$ to ...
Fortunately in situations
arising Hodge theory,  agreement is guaranteed by:

\begin{thm}(\cite[1.3.16,1.3.17, 7.2.5]{D}). 
If $d_i$ strictly preserve $F_{dir}$ 
for $i =0,1...r$, then the first direct and
recursive filtrations  on $E_{r+1}$ coincide. If $r = \infty$ then these
filtrations coincide with the final filtration on $E_\infty$.
\end{thm}

\begin{remark} The above condition for $i=0$ just says that the filtration
on $Gr_W^\dt A$ induced by $F$ is strict.
\end{remark}

\begin{proof}[Proof of theorem \ref{thm:1}]

The argument goes along the lines presented in \cite{D}.
Axiom HC3 and the exactness of $Gr_F$ and $Gr_{\hat F}$ on the
category of $\C$-Hodge structures
implies that $d_1$ strictly preserves $F_{dir}$ and $\bF_{dir}$. The same
goes for $d_0$ by HC2.
 Thus the first direct and recursive filtrations for
$F$ and $\bF$ on $E_2(A, W)$ coincide. Therefore $E_2^{pq}(A,W)$, when
equipped with these filtrations, becomes
 a  pure $\C$-Hodge structure of weight $q$ because it is a subquotient of 
$E_1^{pq}(A,W)$. Consequently 
$d_2 : E_2^{pq}(A,W) \to E_2^{p+2,q-1}(A,W)$ is a
morphism of Hodge  structures 
of different weights and so it must vanish. As zero strictly
preserves everything, the first direct and recursive filtrations on $E_3$
agree, then by similar reasoning   $d_3=0$. Continuing in this
fashion proves (a).

It follows from the previous paragraph that $F$ and $\bF$ induce
 weight $i-k$ $\C$-Hodge structures on $Gr_W^kH^i(A^\dt)$. This implies
 (b). 

The proof of the last  part is based on the fact that given a 
 spectral sequence  of finite dimensional vector
spaces (lying in a translate of the first quadrant), there is an 
inequality 
$$\sum_{p+q=n}\,dim\,E_r^{pq} \ge dim\, H^n(A^\dt)$$
and equality holds if and only 
if $E_r^{pq} = E_\infty^{pq}$ for all $p,q$.
It suffices to prove the opposite inequality for $E_r=E_1(A,F)$. 
Applying the above inequality for a different spectral sequence
gives   
$$\sum_{q}\,dim\,E_2^{q,n-q}(Gr_F^pA, W)\ge  dim\, H^n(Gr_F^pA)
= dim\,E_1^{p,n-p}(A,F).$$
There is a commutative diagram
$$\begin{array}{cccc}
\cdots &E_1^{q,n-q}(Gr_F^pA, W)& \to & E_1^{q+1,n-q}(Gr_F^pA, W) \\
       &||                     &     & ||                        \\
\cdots &H^n(Gr_F^pGr_W^q)      & \to &H^{n+1}(Gr_F^pGr_W^{q+1})  \\
       &||                     &     & ||                        \\
\cdots &Gr_F^pE_1^{q,n-q}(A,W) & \to & Gr_F^pE_1^{q+1,n-q}(A,W)
\end{array} $$
where the top isomorphisms follows from  Zassenhaus' lemma
\cite[1.2.1]{D}, and the bottom from
axiom HC1. As $E_1(A,W)$ carries a  $\C$-Hodge structure and $Gr_F$ is
exact for these, it follows that 
$$ E_2^{q,n-q}(Gr_F^pA, W)\cong Gr_F^p E_2^{q,n-q}(A, W) .$$
By (a),
$$\sum_{q}\, dim\, E_2^{q,n-q}(A,W) = dim\, H^n(A^\dt)$$
and hence
$$\sum_{p}\, dim\, E_1^{p, n-p}(A, F)
 \le \sum_{p,q}\, dim\,Gr_F^p E_2^{q,n-q}(A, W) 
=  dim\, H^n(A^\dt)$$
as required. The argument for $\bF$ is identical.
\end{proof}

\section{Functorial Hodge Structures}\label{sect:fun}

Now I will address the basic problem of constructing
functorial  mixed Hodge structures on the cohomology of a given
 category of spaces and sheaves. 
The obvious strategy is to produce a functorial Hodge complex
for each pair under consideration, but this is usually very
hard to implement. A more useful approach is to 
construct explicit complexes for special pairs, and then try
and approximate the general pairs by the special ones.

Here is the precise set up.
Suppose that $V$ is a category with pairwise fibered products,
a terminal object, and a ``forgetful'' functor (denoted by $|\ldots|$)
to the category $Top$ of 
topological spaces  which preserves these
operations.  Assume furthermore  that for each
 $X\in Ob(V)$ there is a category $\Cf(X)$ equipped with
a forgetful functor (also denoted by $|\ldots|$)
to the category $C^{+}(|X|)$ of bounded below
complexes of sheaves of  Abelian groups  on $|X|$. 
  For each morphism $f:X\to Y$,
there is a pullback functor $f^*:\Cf(Y)\to  \Cf(X)$  compatible with 
$f^*$ for sheaves. One should think of an element of $\Cf(X)$
 as a complex of sheaves 
on $X$
enriched with additional data. 
 $|X|$ will be denoted by $X$
 if there is no danger of confusion.
Let $\Cf$ be the category of pairs $X\in Ob(V), \F\in 
\Cf(X)$, a morphism of pairs is a pair of morphisms $(f:X_1\to X_2, 
f^*\F_2 \to \F_1)$. There is an obvious functor from $\Cf \to V$ and
$\Cf(X)$ is just the fiber over $X\in Ob(V)$.
$\Cf$  will be called  a {\em family of coefficients}.
I will say that $\Cf$ is an {\em Abelian family of coefficients}
if each $\Cf(X)$ is an Abelian category and the functors $f^*$ and
$|\ldots|$ are exact.

Define a Hodge complex on the family of 
coefficients $\Cf$ to be  a functor $A:\Cf \to C_{Hodge}$
such that
$$\begin{array}{ccc}
\Cf& \to & C_{Hodge}\\
\downarrow & & \downarrow\\
  V &\to   & { Top} \\
\end{array}
$$
commutes, and 
a natural quasi-isomorphism between the composite of $A$
with $C_{Hodge}\to C^{+}(-)$ and $|\ldots|$.
In other words, it is just a functorial  choice of a Hodge complex
on $|X|$ naturally quasi-ismorphic to $\F$, for each pair $(X,\F)$ 
with $X\in V$ and $\F\in \Cf (X)$.
 The definition of a
$\C$-Hodge complex on $\Cf$ is similar. 

\begin{lemma} Given a ($\C$-)Hodge complex on $\Cf$, 
$H^i(|X|,|\F|)$ carries a  ($\C$-)mixed structure which is functorial 
in both variables.
\end{lemma}

This lemma is somewhat limited in applicability, because
it's hard to construct Hodge complexes on all of $\Cf$.
In order to formulate a more useful result, let me call 
a full subcategory  $\Delta$ of $V$  {\it dense} provided that:

\begin{enumerate}

\item Given any object $Y\in Ob V$, there is a morphism
$X \to Y$ with $X\in \Delta$ such that the underlying map on spaces
is proper and surjective. Call this  a {\em resolution} of $Y$
in $\Delta$.

\item 
Given a morphism $X \to Y$ and a proper surjective
morphism  $Y' \to Y$,
there is a resolution $X' \to X$ which fits into
 a commutative diagram:
$$
\begin{array}{ccc}
X'& \to  & Y' \\
\downarrow   &&  \downarrow \\
X&      \to & Y 
\end{array}
$$

\end{enumerate}

(It's convenient to phrase the definition
this way,
even though the  second condition implies the first. )
The notion of denseness is extended to 
 coefficients in the most straightforward possible way.
 The preimage of $\Delta$ under $\Cf\to V$ will be called a
subfamily of coefficients; it is called dense if 
$\Delta$
is dense.
The reason for introducing these notions comes from:
\begin{thm}\label{thm:2}
Given a ($\C$-)Hodge complex defined on a dense subfamily of coefficients,
$H^\dt(|X|,|\F|)$ carries  a bifunctorial ($\C$-)mixed Hodge structure for 
all pairs
$(X,\F)\in Ob\, \Cf$. If $\Cf$ is an Abelian family of coefficients then
any short exact sequence of coeffients induces a long exact sequence
of mixed Hodge structures on cohomology.
\end{thm}

An outline of the proof of \ref{thm:2} will be given in 
\ref{sect:desc}.

\begin{ex} Let $V$ be the category of proper complex varieties.
 Let $\Delta$
be the full subcategory of  nonsingular
projective varieties. Density is deduced from either
Hironaka's theorem  \cite{hironaka} (and Chow's lemma),
or the  easier theorem of de Jong \cite{dejong}
on the existence of nonsingular alterations. 
For
$\Cf(X)$, one can take the category  consisting of the constant
sheaf $\Z_X$ and the identity.
 A Hodge complex on $\Delta$ is given by
 \ref{ex:kahler}. This yields the Deligne mixed Hodge structure
on the cohomology of proper varieties (and indeed this is pretty
much Deligne's original construction).
\end{ex}

\begin{ex} Let $V$ and $\Delta$ be as above.
 Let $\Cf(X)$ be the
category pairs $(f:X\to Y, \F)$ where $f$ is a map to a smooth variety
and $\F$ a  polarized complex variation of Hodge structure on $Y$.
 Define $|(f,\F)| = f^*\F$. The $\C$-Hodge
complex of  \ref{ex:vhs}  yields $\C$-mixed Hodge structures
on cohomologies of proper varieties with these coefficients.
\end{ex}

The previous example is essentially due to Lasell \cite{lasell} who
used it in proving his generalization of Nori's Lefschetz theorem
of $\pi_1$. One of the
key steps amounts to verification that if $f:Z\to X$ is a morphism of 
smooth projective varieties (with $X$ and $Y= image(f)$ connected) then
$$ker[H^1(X,\V) \to H^1(Y, \V)] =
ker[H^1(X, \V) \to H^1(Z, f^*\V)]$$
for any polarized complex variation of Hodge structure $\V$ on
$X$.  This 
follows rather formally from the existence and basic properties of 
the $\C$-Hodge structure on these cohomology groups.

\begin{ex}
Let  $V$  be the category of pairs
$(X,Z)$ consisting of a proper complex variety $X$ 
 and a Zariski closed subset $Z$. A  morphism 
$(X_1,Z_1)\to (X_2, Z_2)$ is a morphism of the underlying schemes
such that $Z_1$ is the preimage of $Z_2$. $\Delta$ is the full
subcategory generated by pairs consisting of smooth 
projective
varieties and divisors with normal crossings. Once again density
 follows from  \cite{hironaka} or \cite{dejong}. The Hodge
complex described
in  \ref{ex:log} (with the integral structure described in the
discussion following theorem \ref{thm:1})
yields the Deligne  mixed
Hodge structures on the cohomology of varieties
with constant coefficients.
\end{ex}

\begin{ex} With the same $V$ and $\Delta$ as in the previous
example, let
$\Cf(X,Z)$ be the category of unitary local systems
on $X-Z$. This is an Abelian family of coefficients.
 Then the $\C$-Hodge complex of example
\ref{ex:ulog} gives a $\C$-Hodge structure on the
cohomology of varieties  with coeffients in unitary local
systems. This is compatible with the cohomology long exact
sequence by the theorem, but it
can be also deduced directly from the semisimplicity
of $\Cf(X,Z)$.
\end{ex}

A more elaborate example will be given in 
section \ref{sect:constr}.

\section{Cohomological Descent}\label{sect:desc}

I will outline the proof of theorem \ref{thm:2}.
The name of the game is cohomological descent,
which is a way of computing 
the cohomology of space in terms of the cohomology
of a (simplicial) diagram of spaces lying over it.
It is a generalization
of \v{C}ech's procedure for computing cohomology,
where the diagram in question is just the nerve
of the open cover.
Recall that a strict 
simplicial object in some category $C$ 
is a collection of objects $X_n$ indexed by the natural
numbers, and ``face'' morphisms $\delta_i : X_n\to X_{n-1}$ $i=0,1\ldots n$
such that $\delta_i\delta_j =\delta_{j-1}\delta_i$
for $i < j$. 
I will usually  drop the adjective
``strict'', although this is technically incorrect\footnote{ Strictly
speaking, a simplicial object  has additional  
degeneracy maps $X_{n-1}\to X_{n}$. In all our examples, $C$ will
have finite coproducts, so we can turn a strict simplicial
object $X_\dt$ in a true simplicial object $X_n^+$ by
inductively replacing $X_n$ by $X_n\coprod  X_{n-1}\ldots$, one
$X_{n-1}$ for each degeneracy map... }.
An augmentation  to $X$ is
a collection of morphisms $\epsilon_n: X_n\to X$ which commute with the
face maps. A cosimplicial
object is the dual notion where arrows go the other way, or
equivalently  a simplicial object in the opposite category. Suppose that 
$A^0 \rightrightarrows A^1\ldots $ is a cosimplicial object in an
additive category, then defining $d = \sum_i\, (-1)^i\delta_i$ results
in a complex $A^\dt$.

  Suppose that   $\epsilon_\dt : X_\dt \to X$ is an
augmented simplicial space, i.e. an augmented simplicial
object in the category of topological spaces. A sheaf on $X_\dt$ is a 
collection
of sheaves $F^n$ on $X_n$ and morphisms $\delta_i:\delta_i^*F^{n-1}\to F^n$
satisfying the simplicial identities.  For example,
given an augmentation $\epsilon_\dt : X_\dt \to X$ and a sheaf $F$ on $X$,
the pullback  $\epsilon_\dt^*F$ with $\delta_i^*=id$
is a sheaf on $X_\dt$. 
Suppose that $F$ is a sheaf of Abelian groups.
Then the Godement flasque
resolution $\G^\dt(F^\dt)$ forms a complex of sheaves on $X_\dt$.
Taking direct images yields a cosimplicial object $\epsilon_*(\G^\dt(F^\dt))$
in the category of complexes of sheaves on $X$. This in turn determines
a double complex, and let $s\epsilon_*(\G^\dt(F^\dt))$ denote the total or
associated single complex (I am being a bit cavalier about signs;
the important thing is that they be chosen in a consistent way so as to
yielda complex). The map 
$\epsilon : X_\dt \to X$ satisfies {\em cohomological descent} if
the natural map
$$ F \to s\epsilon_*(\G^\dt(\epsilon^*F))  $$
is quasi-isomorphism for any sheaf $F$ on $X$. 
In less prosaic terms, this says that the natural
map $F\to \R\epsilon_*\epsilon^*F$ is an isomorphism
in $D^+(X)$.

This phenomenon has  already been encountered
earlier  in this article.
Let
$$
\begin{array}{ccc}
E& \stackrel{i}{\hookrightarrow} & \tilde X\\
\downarrow&        & \downarrow \pi \\
S& \stackrel{i}{\hookrightarrow} & X\\
\end{array}
$$
be the diagram given in  proposition \ref{prop:mayervit}. Then
the proof of that proposition actually showed:

\begin{lemma}\label{lemma:mayervit}
$$  
E
\begin{array}{c}\to\\ \to\\  \end{array} 
\tilde X\coprod S \to X
$$
satisfies cohomological descent.
\end{lemma}

The next step is to find a general effective criterion for cohomological
descent. Suppose that $X_\dt$ is a simplicial set, and $x\in X_n$.
If $x_i = \delta_i x$, then these satisfy $\delta_j x_i = \delta_i
x_{j+1}$ for all $i \le j$. Let me call such an $(n+1)$-tuple 
in $X_{n-1}$ an $n-1$ cycle, and let $Z_{n-1}(X_\dt)$ denote the 
set of all of these. Note that $Z_{n-1}(X_\dt)$ depends only on the maps
$\delta_i: X_{n-1}\to X_{n-2}$. This observation will be useful for the 
inductive 
constructions later on. If $X_\dt\to X$ is an augmented simplicial
set, then define $Z_{n-1}(X_\dt\to X)$ by imposing the additional
constraint that the $x_i$ take the same value in $X$.
If $X_\dt \to X$ is an augmented simplicial space  then
 $Z_{n-1}(X_\dt\to X)$ is also a topological space with a continuous map
to $X$. Say that an augmented simplicial space
$X_\dt \to X$ is a  {\it hypercover}\footnote{It would be
 more accurate to call this a hypercover for the Grothendieck topology of
proper maps.}  if $X_0\to X$ and 
the obvious continuous maps $X_n\to Z_{n-1}(X_\dt\to X)$ are  proper and
surjective for each $n$.

\begin{thm} A hypercover satisfies cohomological descent.
\end{thm}

This result will be sufficient for our needs, although it is far from
optimal: it does not include lemma \ref{lemma:mayervit} for
example.
For a more general statement and proofs, see 
\cite[V bis]{sga4} or \cite{gnpp} (this last reference takes a 
slightly different point of view  
into 
which lemma \ref{lemma:mayervit}
fits quite naturally).
Note that
these references use the  stronger form of simplicial object  
having also {\em degeneracy maps}, 
but this will not affect the proofs
since (in the notation of footnote 1) we can pass
from $X_\dt\to X$ to $X_\dt^+\to X$.
Also, in 
order to convert our definition to those in the references, 
$Z_{n-1}(X_\dt)$ is often written as $(cosk\, sk_{n-1} X_\dt)_n.$

Let $V$ be a category satisfying the conditions of the
previous section and let $\Delta$ be a dense subcategory. 
Let $Res$ be the set of resolutions
in the sense of  Section 4
i.e. morphisms with domain
in $\Delta$ such that the underlying continuous map is proper
and surjective.
Given an augmented simplicial object $X_\dt \to X$
in the category $V$,
$Z_{n-1}(X_\dt\to X)$
can be defined as the object that represents the functor
$$T  \mapsto  \{ (f_0,\ldots f_n) \in Hom(T, X_{n-1})^{n+1} | \forall \, 
i \le j,\, \delta_jf_i = \delta_if_{j+1} \}$$
in the category of objects over $X$ (the  representability of this
 functor follows from our assumptions about $V$).
I will treat $Z_{n-1}(X_\dt\to X)$ as an object of
 $V$ by forgetting the map to $X$.
 An augmented simplicial object $X_\dt\to X$ in $V$, will
be called a {\em simplicial resolution } of $X$ if $X_0 \to X$ and each 
$X_n \to  Z_{n-1}(X_\dt\to X)$ lies in $Res$.
A simplicial resolution is a hypercover.

\begin{thm}\label{thm:simpres}
 Every object possesses a simplicial resolution. Given
a hypercover 
$Y_\dt \to Y$ and a morphism $X\to Y$, a
simplicial resolution $X_\dt \to X$ can be chosen so as to fit into a 
commutative diagram
$$
\begin{array}{ccc}
X_\dt&      \to & Y_\dt \\
\downarrow   &&  \downarrow \\
X&      \to & Y 
\end{array}
$$
Any two simplicial resolutions are dominated by a third simplicial
resolution.
\end{thm}

In a nutshell, the proof is  by induction.
The key point is that once $X_i$ and the associated maps are defined 
for all $i < n$, then  $Z_{n-1}(X_\dt\to X)$ exists. And by its very
construction, it carries morphisms to $X_{n-1}$ satisfying the
appropriate face relations.
Now take $X_n$ to be a resolution of $Z_{n-1}(X_\dt\to X)$ (which fits
into a commutative square
$$\begin{array}{ccc}
X_n&  \to& Y_n \\
\downarrow& & \downarrow \\
Z_{n-1}(X_\dt\to X)& \to& Z_{n-1}(Y_\dt\to Y) \\
\end{array}
$$
for the second part).
To start the whole process, set $X_0 $ equal to any resolution of
$X$ (compatible with $Y_0\to Y$). 
For the last statement, apply the previous part with $Y=X$ and
$Y_\dt$ equal to the fiber product of two given simplicial resolutions
of $X$.

Now, I can outline the proof of theorem \ref{thm:2}. To simplify
the discussion, let me dispense with the $\Z$-structure. Suppose
that $(X, \F)\in Ob \Cf$, then choose a simplicial resolution $X_\dt\to X$.
The assumptions, guarantees that there are $\C$-Hodge
complexes $(A^{\dt n}(\F),W,F,\bF)$ quasi-isomorphic to
$\epsilon_n^*\F$ on each $X_n$, and face 
morphisms $\delta_i$ between them. Replace these complexes by 
their Godement resolutions. This results in a differential
graded $\C$-Hodge complex
$$\Gamma(A^{0\dt}(\F))\stackrel{\delta_0-\delta_1}{\longrightarrow}
 \Gamma(A^{1\dt}(\F))\to \ldots $$
The total complex $S^\dt(\F)$ is a $\C$-Hodge complex (example
\ref{ex:total}) which is 
naturally quasi-isomorphic to $\R\Gamma(\F)$. This yields a $\C$-Hodge
structure on $H^\dt(X,\F)$ which depends, apriori,
on the choice of the simplicial resolution $X_\dt\to X$. 
However given a second
resolution $X'_\dt \to X$ dominating $X_\dt$, one gets
a morphism of the corresponding
$\C$-Hodge structures which is an isomorphism
of vector spaces and hence of Hodge structures,
because the forgetful functor to vector spaces is exact
and faithful.
As any two simplicial resolutions are dominated by a third,
this shows independence. Similarly since any morphism can be
lifted to a morphism of simplicial resolutions, this proves
functoriallity. Finally, suppose that $\Cf$ is Abelian, and let
$$ 0\to \F_1\stackrel{f}{\to} \F_2\stackrel{g}{\to} \F_3\to 0$$
be exact. Then there exists  morphisms $S(f):S^\dt(\F_1)\to S^\dt(\F_2)$
and $S(g):S^\dt(\F_2)\to S^\dt(\F_3)$ of Hodge complexes. As 
$S(g)\circ S(f) =0$, there is a well defined map 
$Cone(S^\dt(f))^\dt\to S^\dt(\F_3)$
which  is a 
(nontrilftered) quasi-isomorphism. 
Therefore the $\C$-Hodge
structures on cohomology coincide. The maps in the triangle
$$ S^\dt(\F_1)\to S^\dt(\F_2)\to Cone(f)^\dt\to S^\dt(\F_1)[1]$$
yields maps of Hodge structures in the long exact sequence for 
cohomology.

The following construction is implicit in the above argument. Let
$D_{\C-Hodge}(pt)$ be obtained by localizing the category
of $\C$-Hodge complexes $C_{\C-Hodge}(pt)$ 
with respect to (not necessarily filtered) quasi-isomorphisms
 (see \cite{verdier} for the definition).
Then this  has the structure of a triangulated category with 
translation given  in \ref{ex:trans} and distinguished triangles
represented by
$$A^\dt \stackrel{f}{\longrightarrow}
B^\dt \to Cone(f)^\dt \to A^\dt[1]$$
where $Cone(f)$ is defined in \ref{ex:cone}. There
is a $\Delta$-functor from $D_{\C-Hodge}(pt)$ to the category
of $\C$-Hodge structures which are compatible with the usual
cohomology functor $H^\dt:D^+(\C)\to \C{\textstyle  -Mod}$.
Similar remarks apply to the category of Hodge complexes.
See \cite{beilinson} or \cite{E} for further
details about the  construction. 

\section{ Constant Constructible Sheaves}\label{sect:constr}

Recall that a sheaf of Abelian groups
$\F$ on a complex variety $X$ is
{\em constructible} if
it has finitely generated stalks and 
there exists a partition of $X$ into Zariski locally closed
subsets such that $\F|_U$ is locally constant for each member $U$ of
the partition. I will say that $\F$ is {\em constant constructible}
if each restriction is in fact constant.
As an illustration of the previous ideas, I
want to prove\footnote{This argument was worked at MSRI 
during 1988-1989; my thanks for their support.}:

\begin{thm}\label{thm:constr}
         Let $\F$ be a constant constructible sheaf on a proper 
complex variety $X$, then $H^i(X, \F)$ carries a  mixed Hodge 
structure which is functorial in both variables. A short exact
sequence of sheaves induces a long exact sequence of mixed
Hodge structures on cohomology.
\end{thm}

Let $X$ be a variety over $\C$.
 A {\em stratification} is a finite Boolean algebra $\B$ of subsets of
 $X$ such that the Zariski closure of any element of $\B$ is
also in $\B$. This somewhat unorthodox usage of terminology requires 
some explanation. Fix a stratification $\B$, and let $A(\B)$ be the set
of atoms or nonempty minimal elements. As $\B$ is a finite Boolean
algebra, any nonempty element is a union of atoms. Let me define
 a 
 partial order on $A$, by $Z_1\le Z_2$ if and only if $\bar Z_1
\subseteq  \bar Z_2$. Then, one easily obtains

\begin{lemma} If $Z\in A(\B)$, then $\bar Z = \cup_{Z'\le Z}Z'$ and
$Z = \bar Z - \cup_{Z' < Z}\bar Z'$.
\end{lemma}

Consequently,  $A(\B)$ is a partition of $X$ into Zariski 
locally closed sets, and
the elements of $\B$ are constructible. 
 Conversely, given such a partition, the set of
unions of elements of the partition yield a stratification.

Let $\SV$ be the set of all  pairs $(X,\B)$ consisting of a 
variety with a stratification. This  can be made into 
a category by taking as morphisms $f:(X_1,\B_1) \to (X_2,\B_2)$,
morphisms $f:X_1\to X_2$ of varieties such that $f^{-1}(Z)\in \B_1$
for all $Z\in \B_2$. A morphism corresponding to the identity on the
underlying variety is just a refinement of the stratification. 
A sheaf $\F$ on $X$ with finitely generated stalks will be called
constant constructible
with respect to a stratification $\B$ provided that $\F|_U$ is
constant for each $U\in A(\B)$.

Given a  (Weil) divisor $E$ on 
an irreducible variety $X$, let $\B(E)$ be the smallest stratification 
containing the irreducible components of $E$. If $X$ is smooth and $E$ a 
divisor with normal crossings, then $\B(E)$ will be refered to as a good 
stratification. Let $\Delta$ be the full subcategory of $\SV$
whose objects are disjoint unions of smooth projective varieties with good 
stratifications.

\begin{prop}\label{lemma:stratres}
        $\SV$ satisfies the assumptions of section \ref{sect:fun},
and $\Delta$ is dense.
\end{prop}

\begin{proof} 

We verify that $\SV$ satisfies the axioms of \ref{sect:fun}.
$\SV$ certainly contains a terminal object. Given two
morphisms $f_i:(X_i,\B_i)\to (Y,\B)$ with $i=1,2$, let $P= (X_1 \times_Y 
X_2)_{red}$.
Let $\B_P$ be the Boolean algebra generated by fibered products
$(Z_1 \times_Y Z_2)_{red}$ with $Z_i \in \B_i$. Then $(P,\B_P)$ is the 
fibered 
product of $f_1$ and $f_2$ in $\SV$.

The next step is to show that every object possesses a resolution
with respect to $\Delta$.
If $(X,\B)$ is a stratified variety, let $X'$ be the disjoint 
union of the irreducible components. We can pull $\B$ back to $X'$
to obtain a surjective proper morphism $(X',\B') \to (X,\B)$. Since 
we can work with one component at a time, we may assume that $X$ is 
irreducible. The proof will proceed by induction on the cardinality of
$\B$. If $\B$ has less than or equal to $2$ elements, in other words 
if the stratification is trivial, then the existence of a resolution
follows immediately from \cite{hironaka} or \cite{dejong}. 
In general, let $A$ be the set of atoms of $\B$. Choose an 
element $Z\in A$, which is minimal with respect to $\le$. Then $Z$ is
necessarily closed. Let $\B'\subset \B$ be the Boolean subalgebra
generated by the closures of elements of $A- \{Z\}$. $\B'$ is easily 
seen to be a stratification not containing $Z$. Therefore, by 
induction, there is a resolution $(X_1,\B(E_1) )\to (X,\B')$. Let 
$Z_1$ be
the preimage of $Z$. Then there exists a proper surjective morphism 
$X_2\to X_1$
from a second nonsingular variety such that the preimage $E_2$ of 
$Z_1\cup E_1$ is a divisor with normal crossings [loc. cit.]. 
Then $(X_2,\B(E_2))$ is a resolution of $(X,\B)$.

 It remains to check the last condition for denseness.
Given a morphism $(Y_1,\B_{Y_1})\to
(Y_2,\B_{Y_2})$ and a proper surjective morphism
$(X_2,\B_{X_2})\to (Y_2,\B_{Y_2})$,
let $(X_1,\B_{X_1})$ be a resolution of the fibered product.

\end{proof}

 Given an object $(X,\B)$ of $\SV$, let
$S(X,\B) = \coprod_{U \in A(\B) }\, U$ and let
$J(X, \B):S(X, \B)  \to  X $ denote the canonical map. Set  
$$T=T(X,\B) = J(X,\B)_* \circ J(X,\B)^* :Sh(X) \to Sh(X),$$
where $Sh(X)$ is the category of sheaves on $X$.
From the adjointness of $J_\ast$ and $J^\ast$, there are natural 
transformations:
$$ \eta : I \to T,\, \mu : T^2 \to T$$
satisfying the usual identities of a ``triple" or monad
\cite{maclane}.
Define the natural transformations:
          $$\delta_i = T^i\eta T^{n-i} : T^n \to T^{n+1}$$
 For any sheaf $\F \in Sh(X)$, these maps determine a cosimplicial 
sheaf augmented over $\F$:

        $$ \F  \to  T\F  
\stackrel{\textstyle{\to}}{\to} T^2\F
\ldots $$
and hence also a  complex (augmented over $\F$):
        $$ \F \to  T\F \stackrel{d}{\to} T^2\F \dots$$
where $d = \sum (-1)^i\delta_i$. Either of these will be
denoted by $T^\dt \F$ or $T{(X,\B)}^\dt \F$. 
This can be made  more explicit as follows. Let $j_U:U\to X$ denote
the inclusion for $U\in A(\B)$.
Since $j_V^*j_{U*} = 0$ unless $V \le U$, it follows that:
$$ T^n\F = \bigoplus_{U_0\le \ldots \le U_n}\, j_{U_0*}^*\ldots 
j_{U_n*}^*\F$$
where I abbreviate $j_{U*}j_U^*$ by $j_{U*}^*$. 
The $i$th face map $\delta_i$ is the sum of the
adjunction maps:
$$ 
\ju {0} \ldots \hat \ju {i} \ldots \ju {n} \F \to  \ju {0} \ldots 
\ju {n} \F$$

\begin{lemma}  $T^\dt \F$is quasi-isomorphic to $\F$.
\end{lemma}

\begin{proof}
 Fix $x \in X$. It suffices to construct a functorial retraction 
$h: (T\F)_x \to \F_x$  to the map $\F_x \to (T\F)_x$, because 
then the augmented complex
$ \F_x \to (T^\dt \F)_x$
carries a contracting homotopy defined by 
$$h:(T^{n+1}\F)_x = (TT^n\F)_x \to (T^n\F)_x.$$
See \cite[Appendix 5]{godement} for details. We have
$$ (T\F)_x = \bigoplus_{U\in A(\B)}\, (j_{U*}^*\F)_x$$
If $x \in U$ then $(j_{U*}^*\F)_x = \F_x$, define $h$ to act as
the identity on each
summand $(j_{U*}^*\F)_x$ with $x \in U$, and zero on the other summands. 

\end{proof}

The above constructions are functorial. Given a morphism 
$f: (X, \B) \to (X', \B')$ and a 
sheaf $\F$ on $X'$, there is a morphism of cosimplicial sheaves:
$$ f^*T(X',\B')^\dt \F \to T(X,\B)^\dt f^*\F.$$

\begin{proof}[Proof of theorem \ref{thm:constr}]
Let $\mathcal{PSV}$ be the full subcategory of $\SV$ where
the underlying variety is proper. The intersection
$P\Delta= \Delta\cap Ob\,\mathcal{PSV}$ is dense in $\mathcal{PSV}$.
Let $\Cf(X,\B)$ denote the category of constant constructible
sheaves on $(X,\B)\in Ob\,\mathcal{PSV}$. This determines an
Abelian  family of coefficients.

Let $(X,\B)$ be an object of $P\Delta$, and  
let $\F$ be a constant constructible  sheaf on it. 
There is no
loss of generality in assuming that $X$ is connected, so that
$\B$ is the Boolean algebra $\B(E)$ generated by the components of a
divisor with normal crossings. In this case, the components
of $T^\dt(\F)$ are 
fairly easy to describe: $j_{\bar U_0}^*\ju{1} \ldots \ju {n} \F = 
j_{\bar U_0}^*\ju {n} \F$ is a constant sheaf.  Let 
$$A^{n\dt}(\F) = 
\bigoplus_{U_0\le \ldots \le U_n}\, 
j_{\bar U_0*}(\E_{\bar U_0}^\dt\otimes j_{\bar U_0}^*\ju {n} \F 
,\ldots)$$
be the associated Hodge complex. The differentials of $T^\dt(\F)$
extend to give a differential graded Hodge complex $A^{\dt\dt}(\F)$. 
Then $\F\mapsto tot^\dt(A(\F))$
is a Hodge complex on $\Cf|_{P\Delta}$. The theorem now follows
from theorem \ref{thm:2}.

\end{proof}

There are a couple of variations on this construction.
Call a sheaf $\F$ of $\C$-vector spaces {\em unitary constructible}
if it has finite dimension stalks and there is stratification
such that $\F|_U$ is isomorphic to a unitary local system for each
$U\in A(\B)$.  In this case, the sheaves
$j_{\bar U_0}^*\ju{1} \ldots \ju {n} \F$ are all unitary local
systems, let $(V_{(U_0,\ldots U_n)},\del)$ denote the extension
of the associated flat vector bundle to $\bar U_0$ 
(see example \ref{ex:ulog}). Then
 a $\C$-Hodge complex:
$$A^{n\dt}(\F) = 
\bigoplus_{U_0\le \ldots \le U_n}\, 
j_{\bar U_0*}
(W_0 A_{\bar U_0}^\dt(log(\bar U_0-U_0)) \otimes V_{(U_0,\ldots U_n)
} ,\ldots)$$
can  be constructed on $P\Delta$,
where $A_{\bar U_0}^\dt(log(\bar U_0-U_0)) \otimes V_{(U_0,\ldots
  U_n)}$ is described  \ref{ex:ulog}. This  produces 
a natural $\C$-Hodge structure 
on cohomology with coefficients
in   $\F$.
Similar arguments yield ($\C$-)mixed Hodge
structures on the cohomology of arbitrary varieties with constant
(or unitary) constructible coefficients. The 
bookkeeping  is a bit more
involved in this case. $\SV$ would be replaced by
the category of triples $(X,\B, \bar X)$ where $\B$ is a
stratification of $X$ and $\bar X$ a compactification.  The objects
of $\Delta$  would be triples where $\bar X$ is disjoint union of
smooth projective
varieties, $\B$ a good stratification of the form $\B(E)$
and  $(\bar X - X)\cup E$ a divisor with normal crossings.

I want to conclude 
with a few words about how all of this ought to relate
to Saito's work on mixed Hodge modules \cite{saito1, saito1.5, saito}.
I will treat mixed Hodge module theory entirely as
black box, and just summarize some of the
main properties. For a gentle introduction to the subject
see \cite{brylinski-zucker}. To every 
variety  $X$ over $ \C$,
 Saito constructs the Abelian category
$MHM(X)$ of mixed Hodge modules. For a point, $MHM(pt)$ is equivalent
to the category of $\Q$-mixed Hodge structures $H$ which are polarizable
in the sense that $Gr^W_\dt H$ admits polarizations. 
There is an
additive  functor $rat: MHM(X)\to D^b_{c}(X,\Q)$
to the derived category of sheaves of $\Q_X$-modules with bounded
constructible cohomology; in fact, $rat$ is an exact
and faithful functor to  the Abelian
subcategory of perverse sheaves \cite{bbd}.
When $X$ is smooth, the category
of perverse sheaves over $\C$ is exactly the image of 
the category of regular holonomic $D_X$-modules under
the De Rham functor $DR$ (by the Riemann-Hilbert
correspondence).
In this case, $MHM(X)$ is constructed as a subcategory of 
the category whose objects are pairs consisting of a filtered regular
holonomic $D_X$-module $(M, F)$ and a filtered
perverse sheaf over $\Q$ $(K, W)$ such that $DR(M)\cong K\otimes \C$.
 The  functor $rat$ extends to a morphism of
triangulated categories denoted by the same
symbol $rat:D^bMHM(X)\to D^b_{c}(X,\Q)$.
The basic operations of sheaf theory extend to  $D^bMHM$.
In particular, given a morphism $f:X\to Y$, there exists
morphisms of triangulated categories 
$\R f_*: D^bMHM(X)\to D^bMHM(Y)$
and $\LL f^*: D^bMHM(Y)\to D^bMHM(X)$
which are compatible with the corresponding morphisms of
$D^b_c$  (I have added
the ``$\R$''  and ``$\LL$'' 
for consistency; Saito leaves these out
by convention).
These facts together imply that $H^i(X, rat(M))$ 
is assigned a mixed Hodge
structure for each $M\in D^bMHM(X)$.

Given a constant constructible sheaf $\F$ on $X$, I expect that there
should be a natural choice of  
$M\in D^bMHM(X)$ with $rat(M) =\F\otimes \Q$
such that $H^i(X,M)$ gives the same $\Q$-mixed Hodge structure as above.
When $\F = A_X$ is constant, the obvious choice
for $M$ is $\LL f^*(A\otimes \Q(0))$, and
when $X$ is smooth the  mixed Hodge structures on
$H^i(X,A)$  do coincide by \cite[3.3]{saito1.5}.
 I  have a candidate
for  $M$ for general $\F$, but I don't yet have a proof that it works.
The main problem is that the derived category is the wrong
place for doing descent,  a more ``rigid'' description
of $M$ is needed first. Mixed Hodge modules
do not work on simplicial schemes for these reasons; 
some of these issues
are touched upon in a recent preprint of Saito \cite{saito-new},
but it appears that more work is needed to resolve them.

\end{document}